# On the statistical distribution of prime numbers, a view from where the distribution of prime numbers is not erratic.

**Sandor Kristyan**
**Hungarian Academy of Sciences, Research Center of Natural Sciences**
1117 Budapest, Magyar tudosok krt. 2, Hungary, kristyan.sandor@ttk.mta.hu

## Abstract

Currently there is no known efficient formula for primes. Besides that, prime numbers have great importance in e.g., information technology such as public-key cryptography, and their position and possible or impossible functional generation among the natural numbers is an ancient dilemma. The properties of the functions $2ab+a+b$ in the domain of natural numbers are introduced, analyzed, and exhibited to illustrate how these single out all the prime numbers from the full set of odd numbers. The characterization of odd primes vs. odd non-primes can be done with $2ab+a+b$ among the odd natural numbers as an analogue to the other, well known type of fundamental characterization for irrational and rational numbers among the real numbers. The prime number theorem, twin primes and erratic nature of primes, are also commented upon with respect to selection, as well as with the Fermat and Euler numbers as examples.

## Keywords

prime number generator, prime number theorem, twin primes, erratic nature of primes

## 1.Introduction

*In medias res*: The $m= 4ab+2(a+b)+1$ jumps over all odd primes (except 2) and makes multiple hits on all odd composite (non-prime) numbers for $a,b$=1,2,…, while if zero is also allowed for $a$ and $b$, it generates all odd numbers – this will be introduced and analyzed next. In the literature it is stated that, there is no known useful formula that sets apart all the prime numbers (a natural number greater than 1 that has no positive divisors other than 1 and itself) from composites (non-primes, i.e., a product of two other numbers both greater than 1). Primes are used in several routines in information technology, such as public-key cryptography which makes use of properties such as the difficulty of factoring large numbers into their prime factors [footnote 1]. Legendre has shown us that there is no rational algebraic function that always gives primes. In 1752 Goldbach illustrated that no polynomial with integer coefficients can give a prime for all integer values [1-2]. The best-known polynomial that generates primes only, is the famous $n^2+n+41$ accredited to Euler in 1772 [1] which gives distinct primes for the 40 consecutive integers $n$= 0 to 39, and besides it soon terminates in relation to the fact that there are infinite primes; it even jumps over primes, e.g., for $n$= 3 and 4 it yields 53 and 61 - jumping over prime 59. There are some low-order polynomials that only generate primes (a few dozen) for the first few non-negative values [3, footnotes 2-3]. There are formulas generating the prime numbers, exactly and without exception, but there is no such formula which is known to be efficiently computable. Two of the most famous and earliest are the Fermat ($2^r + 1$ with $r\equiv 2^q$ for $q$=0-4) and Mersenne ($2^p - 1$, where $p$ is assumed prime) numbers [4-5, footnotes 4-5]. The problem with these polynomials and formulas are, that they do not predict all primes, and even worse, they predict only a few primes - general reviews and textbooks on the basic properties of primes can be found in refs [6-11]. Nowadays, in the century of computation and internet data, prime numbers up to 1 trillion is easily available from the web, see C.K. Caldwell (research or data) home pages, etc. [footnotes 6-7]. As a characterization, Wilson's theorem should also be mentioned which states that: A natural number $n > 1$ is a prime number $\Leftrightarrow (n\text{-}1)! \equiv -1 \pmod{n}$.

In spite of these facts, there are still two simple generators predicting all the primes: the trivial (and useless at first glance) first order polynomial generating the set $S_2 \equiv \{2k+1 | k=0,1,2,3,\ldots\}$ of all odd primes and all odd non-primes, i.e., all the odd numbers, as well as the non-trivial, relatively famous and faster [12] polynomial generating the set $S_6 \equiv \{6n \pm 1 | n=1,2,3,\ldots\}$ of all odd primes (except 3) and not all but i.e., all except $6n+3$, see Lemma 1 below, infinite odd non-primes. (Indices 2 and 6 on S refer to the coefficients.) Furthermore, there is only one even prime, the 2 ("queen" of primes) completing the set of primes ("kingdom" of number theory or aritmetics). Below, for simplicity, we mean "odd numbers" the "odd positive numbers" or "odd natural numbers", if "odd negative numbers" are discussed we will emphasize that, also if zero is included or not. We use the term odd primes emphasizing their odd nature, because our discussion below is heavily based on the odd number generator $2k+1$. A minor technical consequence is, that for odd primes $\{3,5,7,11,\ldots\}$ the order numbers are $\{1,2,3,4,\ldots\}$, respectively i.e., the set of natural numbers, that we use here, although for the full set of primes $\{2,3,5,7,11,\ldots\}$ this set of order numbers shifts to the left by one, and not prime 3 but prime 2 is the first.

In fact, generators in $S_2$ and $S_6$ are polynomials wherein the variables $k$ and $n$ are order numbers as well, and to top it all, the simplest ones - first orders. In $S_2$ the $2k-1$ also has a good function, dropping the $k=0$ from the domain; as well as this we consider only positive integers for order numbers, although the negatives will also be commented upon. Below, we introduce simple polynomials - indicated in the title - which, select all the odd non-primes from the sets in $S_2$ and $S_6$ generating all the odd prime numbers as residue. From the set, $S_2$ the $k=2ab+a+b$ selects out the odd non-primes with $a=1,2,\ldots\infty$ and $b=1,2,\ldots a$, or a reverse domain definition in respect to the role of $a$ and $b$ having this function symmetric for $a$ and $b$, and as an alternative, from the set $S_6$ the $n= 6k_1k_2 \pm k_1 \pm k_2$ select out the odd non-primes in a similar manner. Notice that $S_2$ is "primarily called" the odd number generator, while $S_6$ is "primarily called" the odd prime generator – had $S_6$ not generated odd non-primes either, the functional generation of primes would be a simple elementary task, ending the ancient, furious search for a prime number generator throughout the history of number theory. The dream is to find a simple form like the Euler formula above if it exists, because until now, only the ancient, unrefined method exists, such as diagnosing every natural number one by one to see if it is a prime or not – of course, there are many useful theorems and algorithms which speed things up, but it is a very still need laborious process. The latter provides the basic part of the competition over centuries to find the largest (odd) prime in the year published.

## 2.1 The prime selector function *2ab+a+b* for all the odd numbers

The way to this selector function is very simple. The "fundamental theorem of arithmetic" states that all natural numbers (1,2,3,…) can be written exactly in one way, in the form $p_1^{n1}$ $p_2^{n2}$... $p_L^{nL}$, with natural numbers $L$, $n_1$,...$n_L$ >0 and primes $p_1 < p_2 < . . . < p_L$, knowing the contiguous theorem such as, there are infinitely many primes, more precisely, there is one even prime, the 2, and the rest are the infinite set of odd primes in the form of $2k+1$, where $k$ is from the set of certain but not all natural numbers, so, the partition between even prime(s) and odd primes is extremely unbalanced - known even from elementary maths textbooks. Considering only odd numbers, the prime 2 cannot be among $p_i$'s, and for simplicity, keeping the same order, the partitioning of this product can be written as product of two, $m_1m_2$, such as $m_1=p_1^{n1}$ $p_2^{n2}$… $p_Q^{nQ}$ and $m_2=p_{Q+1}^{nQ+1}$ ...$p_L^{nL}$, with the important constrain that both $m_1$, $m_2$ >1, as well as this, these $m_i$'s can be written in the form $m_1=2a+1$ and $m_2=2b+1$ with proper natural numbers $a$ and $b$. We know that the product of two odd numbers in the form $2k+1$ preserves this form, and also that $m_1$ and $m_2$ are odd primes or odd non-primes; e.g. for 2275= $5^2$*7*13 the $(m_1,m_2)$= (5,5*7*13)=(2*2+1,2*227+1) $\Rightarrow$ $(a,b)$=(2,227) or $(m_1,m_2)$= ($5^2$,7*13), etc. i.e., the partitioning is not unique, but that causes no problem in the derivation below. In this way, any odd number can be partitioned as $m=2k+1= (2a+1)(2b+1)= 2(2ab+a+b)+1$, and since the constraint $m_1$, $m_2$ >1, m cannot be prime for $a,b$>0 this latter ensures that $m_1,m_2$>1, and finally we have reached the selector function from the right most side as

$$k= 2ab+a+b \qquad \text{(Eq.1)}$$

for generating or selecting the odd non-primes

$$m= 2(2ab+a+b)+1 = 4ab +2(a+b) + 1. \qquad \text{(Eq.2)}$$

**Theorem 1:**

All odd non-primes can be selected from the set $S_2$ with the help of Eq.2:

$$\{2(2ab+a+b)+1|\ a=1,2,\ldots\infty \text{ and } b=1,2,\ldots a\} = \text{all the odd non-primes}, \qquad \text{(Eq.3)}$$

and those that are jumped over, those are all the odd primes.

For scan, let us call full scan for $(a,b)$ the $a,b$= 1,2,3,…,∞, and restricted scan for $(a,b)$ the $a$= 1,2,3,…,∞, $b$= 1,2,3,…,$a$. We note that interchanging $a$ and $b$ can also be defined as "alternative", also in Eq.3, the restricted scan is used. Besides that, the reversed role of variables $a$ and $b$ tells us the same in Eq.3, being $S_2$ and $S_6$ symmetric in respect to interchange of $a$ and $b$, also, in spite of this restriction for $b$ in relation to $a$, $b$ goes up to a only, because $(a,b)$=(4,2) and (2,4) would yield the same $k$), Eqs.1-2 can still generate the same numbers with different non-symmetric values, e.g. $(a,b)$=(4,2) and (7,1) yield the same $k$=22 and $m=2k+1$= 2*22+1=45 (=5*9) non-prime, but there are many $k$'s generated only once. Full scan can also be used in Eq.3, but in that case the same odd non-primes are generated even more (doubly). The diagonal elements ($a=b$) of Eq.2 are $m= 4a^2+4a+1= (2a+1)^2$, the square odd numbers which are odd non-primes in Eq.3 (see $m$ in Table 1), as well as this, arbitrary large odd non-primes can be generated simply and assertively via Eq.2, for example, with $a$=123456789 and $b$=987654321, the 18 digits $m$=487730526672763297, which is "not a big deal" in the problem of prime generation, since it can be done by multiplying any two big numbers. The other side of the latter, i.e., generating very large (odd) primes by Eqs.2-3 is, in principle, possible because there are infinite primes, but it is not easy and neither is it addressed here i.e., finding a very large $k$ which is jumped over. If $a$=0 and/or $b$=0 was allowed in Eq.3, the set would include all the odd numbers (prime and non-prime) as $k=b$ and $m=2b+1$, expanding the full selected set in Eq.3 back to the full un-selected set $S_2$, - the reason is that $m_1$ and/or $m_2$ become(s) 1 above, violating the $m_1$ and $m_2$ >1 requirement.)

**Demonstration of Theorem 1:**

Table 1 shows some initial values generated by Eqs.1-2, and only the lower diagonal is exhibited for this symmetric matrix; the multiplicity or distribution how the values of $k$ is hit by scanning ($a,b$) is exhibited in Table 2 along with a larger scale graphical representation on Figures 1-2 how the odd primes are selected in a given interval, particularly in [1,100] and [901,1000] as a spectrum, (see also Appendix 1). Going back to the matrix representation (diagonal, $a=b$, and off-diagonal, $a\neq b$, elements m in Table 1), those are the odd primes which are "jumped over" or the "holes" in that matrix.

**Extension of Theorem 1:**

What if in Eq.1 $a,b$= -1,-2,-3,…-∞? Starting again from $m=(2a+1)(2b+1)$, armed with the knowledge that positive and negative odd numbers preserve the $2k+1$ form under multiplication, the $a$=-1 (or $b$=-1) and allowing $b$ (or $a$) to run, Eq.1 becomes $k= -b-1$ (or $-a-1$) and scans all the order numbers (all natural numbers); so, Eq.2 scans all positive odd numbers, as expected since with $a$=-1 the $2a+1$= -1, so it does not force two non-unit numbers in the product for $m$. However, for $a,b$= -2,-3,…-∞ Eq.1 scans the same as Table 1, and the order number of odd primes are jumped over again, so the positive odd non-primes are generated. The matrix in this case in Table 1 shifts diagonally, i.e., the 4 will be in (row, column) position ($a,b$)= (-2,-2) and not at (1,1) as in the case of $a,b$= 1,2,3,…∞ (full scan) before, but the symmetric property also holds for scanning $a$ and $b$ as in Eq.3 and, as a consequence, negative sign analogue of restricted scan is enough. (Notice that $a,b$>0 and $a,b$<0 both yield $k$>0 via sign cancellation in $ab$.); not surprisingly, the $a,-b$ or $-a,b$= 1,2,3,…∞ cases yield the same selections via Eq.1 in the same manner, but generate negative odd numbers. If one does not like negative parameters, the $a:=-a$ and/or $b:=-b$ yield the alternative and equivalent expressions for Eq.1 with respect to selector ability as $k=2ab\pm a\pm b$, but one should be careful how to run ($a,b$), because a bad combination, e.g., $k=2ab+a-b$ with $a,b$=1,2,… yield the useless mixed set of odd primes and non-primes, while $a,-b$=1,2,… provides the same as Eq.3. Furthermore, the entire discussion in this section can only be had if one starts from the $2k-1$ form for odd numbers, instead of the equivalent $2k+1$ form.

### 2.2 Image of odd non-prime generator $m= 2(2ab+a+b)+1$ for total selection, particular example for Fermat and Euler numbers using selector function $2ab+a+b$

Eq.1 selects all the order numbers ($k$) from the set of odd numbers belonging to non-primes ($m$ in Eq.2), via the set and domain ($a,b$) in Eq.3. We emphasize that while Eq.1 is a relatively simple algebraic form, it totally or fully selects for $S_2$. The function $m$ in Eqs.2-3 generates all odd non-primes, and jumps over all odd primes, so the complement of its image (complement Im ($m$)) in the set of odd numbers are the odd primes, and these numbers or subsets of these are targeted to be described as Fermat, Mersenne, numbers etc.. Preliminary analysis has shown that although Eq.2 has a relatively simple functional form, it is a very difficult task to describe those numbers (odd primes) with simple functions - the $m$ jumps over, (see Introduction). To be more precise, the old dilemma of the existence and algebraic form of the prime number generator is not yet stated or solved, but the complement set in Eq.3 with respect to odd numbers via the simple Eq.2 is established. Particularly, to find expression what Im(*m*) in Eq.2 does not hit is not solved by the author or anybody else, since to the knowledge of the author the function in Eq.2 has not been analyzed in this respect, so it may show a useful way to search further. In other words, the search for a prime generator can be achieved by accomplishing algebraically the complement image of function *m* in Eq.2 with respect to odd numbers.

**Application of Theorem 1 for Fermat numbers:**

For the subset of these *m*, one can recall the old example of Fermat numbers which are ($F(q) \equiv 2^r + 1$ with $r \equiv 2^q$ for $q$=0-4) $2^1+1=3$, $2^2+1=5$, $2^4+1=17$, $2^8+1=257$, $2^{16}+1=65537$ mentioned in the Introduction, and are primes (belonging to the complement image of function *m* in Eq.3), but there are no other known Fermat primes with $q > 4$, i.e., those are odd non-primes belonging to the image of function *m* in Eq.3. One can see the general fault of current prime generator functions in the case of Fermat numbers too; not only does it stop very soon, namely at $q$=5, but it also jumps over odd primes when $q$<5, e.g. the prime 7, etc..) The hypothesis for Fermat numbers such as F($q$) is prime for $q$=0-4 as simple calculation proves, but non-primes for $q$>4 is not an easy task to prove and is not addressed here. However, with Eq.2 the $m \equiv 2(2ab+a+b)+1 = 2^r + 1$ with $r \equiv 2^q$ can be analyzed here in a short note, which rewrites this hypothesis as always $\exists$ ($a,b$) pairs for $q$>4 to satisfy this equation, but no ($a,b$) exists for $q$=0-4. This equation can be simplified as

$$k \equiv 2ab+a+b = 2^{r-1} \quad \text{with } r \equiv 2^q \ . \qquad \text{(Eq.4)}$$

Right hand side of Eq.4 are 1, 2, 8, 128, 32768 for $q$=0, 1, 2, 3, 4, respectively. Table 1 shows that $k$= 1, 2, 8 are jumped over, so are the 128 and 32768, but they are not exhibited, because there are no ($a,b$) yielding these numbers, to tell us that these F($q$) are primes, but of course, this check can be done in other regular ways. For $q$>4, the hypothesis says that $\exists$($a,b$) always pairs to satisfy Eq.4, telling us that these F($q$) are non-primes; e.g., for $q$=5, the ($a,b$)= (3350208,320) calculation search or scan (of which the complexity is rapidly increasing with $q$) yields $k$= 2147483648 = $2^{31}$ for the left in Eq.4, which is the right side also for $q$=5; and indeed, $m$= 2(2*3350208*320+3350208+320)+1 = 4294967297 = 641*6700417, a non-prime.

**Remark on Theorem 1 in relation to Euler primes:**

The Euler primes (see Introduction) can be exemplified in just the same way as for Fermat numbers, but for the sake of brevity, we will not make a similar example calculation but instead, we will make a note as a hypothesis. The Euler formula $n^2+n+41$ generates some, (40 pieces of small) odd primes, while Eq.2 generates all the (infinite) odd non-primes, so 40 elements of the above mentioned complement set of Im(*m*) can be described by the Euler formula. We draw attention to that search, scan or calculation for the complement of Im(*m*) in Eqs.1-2, for odd primes, as demonstrated in Figures 1-2 for K=0 level is more powerful then e.g, the Euler formula for (40) odd primes in view of the existence of generator function, because Eqs.1-2 selects all the infinite ones. The simple reformulation of Eq.2, $m= 4ab+2(a+b)+1$, for all odd non-primes is similar to the Euler formula in that both are second order polynomials, and dissimilar in the manner that they predict complement set elements (odd non-primes vs. odd primes), Eq.2 is a two variable polynomial in contrast to the one variable Euler polynomial, and, last but not least, Eq.2 makes total selection for infinite elements while the Euler formula predicts only 40 at the beginning.

**Weak conjecture after Theorem 1:**

In relation to Goldbach's statement in the Introduction, if one seeks polynomial description for odd primes, it must be at least two dimensional (or variable) and at least second order polynomial, since the set in Eq.3 (all odd non-primes) cannot be described with one dimensional polynomial; so with this in mind, there is little hope that one variable polynomial can describe the complement set, (obviously, e.g., $S_2$ or $S_6$ describe it with the first order one dimensional), but without selection. The simple and powerful form in Eq.3 may indicate that two dimensional polynomials may exist generating many more odd primes than the Euler formula, or it would be wonderful if it was a low order one for generating all odd primes (see some more notes in Appendix 2). Complex numbers (containing the powerful imaginary unit $\sqrt{-1}$) may also need to be involved, however, here we will not address this question any further.

**Corollaries of Theorem 1 in relation to emblematic theorems on primes:**

The theorems such as "There are infinite primes…", "The density of primes decreases…", "For any natural number $n > 1$, there exists an (odd) prime number $p$ such that $n < p < 2n$ (Chebyshev)", "Every even integer $2n$ greater than 2 can be expressed as the (not necessarily unique) sum of two primes (Goldbach's conjecture)", etc. can be composed or stated for the image set of Eqs.1-2 as "There are infinite natural numbers that $k$ and $m$ jump over" because those are the order numbers and values of odd primes, respectively. "The density of the $k$ and $m$ values increases as $a$ and $b$ increase", "For any natural number $n > 1$, there exists an integer $p$ (not necessarily only one) such that $n < p= 4ab +2(a+b) + 1 < 2n$ has no integer solution for $a,b>0$, except when one and only one in a pair $\{a,b\}$ is zero; this $p$ is an odd prime (otherwise, if $a$ or $b$=0, but not both, the $p$ scans all odd numbers in $(n,2n)$)", etc.

**2.3 Characterization of odd non-primes vs. odd primes with the selector function**

The definition of prime numbers can be found in the beginning of the Introduction. The $2k+1$ in $S_2$ generates all positive odd numbers, of which a non-empty infinite subset contains all positive odd primes. Besides $2k$ generates the zero ($k$=0) and all positive even numbers including the only one positive prime 2 ($k$=1), the $2k+1$ generates the unit element ($k$=0) also, a handicapped value $k$= ½ generates the only one positive prime 2. The $m= 2(2ab+a+b)+1$ in Eq.2 generates all the positive odd non-primes, a countable subset, with the $(a,b)$ values in Eq.3. that are jumped over, those are all the positive odd primes, also a countable subset – allowing a characterization as follows:

**Theorem 2:**

Eq.3 partitions the set of natural numbers $\{a|a=1,2,3,\ldots\infty\}$ to odd primes and the remainder.

**Remark on Theorem 2:**

The properties of such an $m$ in Eq.2 have been analyzed and demonstrated above. In this respect the known partition of the set of real numbers should be recalled as analogy, in which set the rational numbers (countable subset) are defined or generated by $a/b$ with $a,b= \pm1,\pm2,\ldots$ ($a$=0 is also allowed), those that are jumped over are the irrational numbers (uncountable subset), based on theorems that $\sqrt{2}$, e=2,718… or π=3.1415… etc. cannot be written in the form of $a/b$. Even more is true in this analogy, the $a/b$ (selector function of all rational numbers from irrational ones in the set of real numbers) and $m= 2(2ab+a+b)+1$ in Eq.2 (selector function of all odd non-primes from odd primes in the set of odd numbers) are both invariant on $a \leftrightarrow b$ interchange, and of course $a/b$ switch to a reciprocal value and $b$=0 is allowed instead but, remains rational in nature, while in Eqs.1-3 even the expressions remain the same, especially if $b$ is not restricted for practical purposes as in Eq.3, as well as $a/b$; also Eq.2 can provide the same values multiple times, e.g., 1/11=2/22 for rational numbers, or see the two appearances of number $k$=22 in Table 1 or the K>1 peaks in Figures 1-2.

It is interesting to note that a simple rearrangement of Eq.2 yields $m= [2(a+b)+1] + 4ab$, showing that full scan in the square bracket generates all the odd numbers multiple times, e.g., $(a,b)$=(1,2) and (2,1) or (2,8) and (3,7) generate the same 7 or 21, respectively with the exception of 1 and 3 (via $a+b$=0 and 1, respectively) and the superposition or additive term $4ab$ generates, selects or reduces this set to all the odd non-primes for $m$. An elementary real function property is recalled to compare with the selecting role of the second order term $4ab$ as an analogy; e.g. on $x \varepsilon (-\infty,\infty)$ the image of $x+0.25$ is $(-\infty,\infty)$, while $x^2+x+0.25$ selects the negative real numbers out, or reduce this image set to $[0,\infty)$ by the selecting property of the second order term $x^2$. An immediate thought is that Eq.2 (rearranged in the second order form) $m= 4ab +2(a+b) +1$ always allows us to write the product form $(2a+1)(2b+1)$ for $a,b$ > 0, but Eq.1 jumps the order number of odd primes; similarly as the second order real coefficient algebraic equation $y= x^2 +bx +c$ allows the product decomposition $(x-x_1)(x-x_2)$, but with the known selection as if discriminant $b^2 -4c \geq 0$ the $x_i$ remain in the real set otherwise it is in a complex set, i.e., product decomposition is not possible in a real set. A similar thing happens in Eqs.1-2 wherein product decomposition $m=m_1m_2$ is not possible for $m_1,m_2$>1 at certain $k$'s, and to search the detailed analogy further is worth doing. The important thing in this analogy is that Eqs.1-2 with its jump property over odd primes is similar to the second order real coefficient algebraic equation (starting at $x^2+1=0$) how it encounters real vs. complex numbers, that is to say, how Eqs.1-2 encounter odd non-primes vs. odd primes.

### 3.1 The prime selector function $6k_1k_2 \pm k_1 \pm k_2$ for odd numbers of the form $6n \pm 1$

First of all, we introduce another type of selector function:

**Lemma 1:**

Reducing the set from $S_2$ to $S_6$ is done with a certain selector function, $6n+3$, selecting out infinite but not all odd non-primes.

Furthermore, proving the $6n \pm 1$ rule for the odd prime generator in $S_6$ is very easy, and we show simpler and more convincing proof here than in ref.[12].

**Proof of Lemma 1:**

Calling $6n-1$ the lower value and $6n+1$ the upper value, any $n$ generates two neighbor odd numbers by jumping over an odd non-prime number in the interval between the upper ($6n+1$) and $n+1$ generated lower ($6(n+1)-1= 6n+5$) values, i.e., in the interval $[6n+1, 6n+5]$, that is in the odd integer set $\{6n+1, 6n+3, 6n+5\}$. Actually, only the odd numbers in the middle, the $6n+3$ are jumped over by $6n \pm 1$ vs. $2k+1$ which are odd non-primes, since $6n+3= 3*(2n+1)$, i.e., the 3 is always the divisor, while its left and right odd neighbors are of the form $6n \pm 1$ odd primes or odd non-primes. In fact, this is the easiest and simplest way to increase the power of the generator function $2k+1$ to $6n \pm 1$ to generate odd primes by getting rid of unwanted $6n+3$ odd non-primes – although it is still a long way from being a full selection. The statement that: “All primes except 2 and 3 is of the form *6n±1*, where *n* is a natural number" is a typical number theory statement in the sense that, it looks extremely difficult to prove at first glance, but in fact it is very easy to prove – one can even find naive computer checks for this form on the internet - checking primes up to a billion, in this respect.

**Remark on Lemma 1:**

Note that *6n+3* is also a certain selector function reducing the set $S_2$ to $S_6$, but it still leaves odd non-primes in the set, while selector function Eqs.1-3 single out all odd non-primes from set $S_2$. Note also, that the *6n±1* rule selects out from *2k+1* the numbers which have the smallest prime divisor 3, similarly, but much more complex less elegant and less compact forms can be derived if primes 3 and/or 5 divisors are singled out in the same manner, etc.. There is another generator known in literature - the *4n±1*, stating by indirect proof [12] that it generates infinite primes; however, scanning via *n*, one realizes that it generates the same set of all odd numbers such as *2k+1*, i.e., falls into set $S_2$; It is highly relevant [12] that both signs in *4n±1* generate infinite odd primes just as in *6n±1*.

**An application of Lemma 1:**

Comparing $S_6$ to $S_2$, *6n±1* selects faster, furthermore, it provides a useful check for any individual odd prime because to decide if the 6 is a divisor is a simple task. that is, half of the number is a multiple of 3, and the latter holds if the sum of digits is a multiple of that number i.e.,3. If the m in $6|(m+1)$ or $6|(m-1)$ holds for a number $\Rightarrow$ *m* has the form *6n±1* $\Rightarrow$ it can be an odd prime. If 6 is not a divisor of *m+1*, nor of *m-1* $\Rightarrow$ *m* is not in the form *6n±1* $\Rightarrow$ but it is an odd non-prime. For example, $m=17 \Rightarrow 6|(17+1) \Rightarrow$ it can be prime, actually, it is, or if $m=15 \Rightarrow$ 6 is not a divisor of 15+1, nor of 15-1 $\Rightarrow$ not a prime, indeed, it is not (15=3*5); but at this stage we shall not go into details on how a number is factorized, neither on how the difficulty increases for large numbers.

Eq.2 and the generator in $S_6$ yield that if integers *n, a* and *b* satisfy

$$2(2ab+a+b)+1= 6n\pm1, \qquad \text{(Eq.5)}$$

then the right hand side must be odd non-prime giving two immediate conditions because the left hand side itself is always odd non-prime, and the right hand side itself is generally odd prime or non-prime so, the left hand side limits the two choices to one:

$$n=(2ab+a+b)/3 \qquad \text{yields all odd non-primes of form } 6n+1 \qquad \text{(Eq.6)}$$

$$n=(2ab+a+b+1)/3 \qquad \text{yields all odd non-primes of form } 6n-1 \qquad \text{(Eq.7)}$$

for integer *n* values. For scan, let us call $(k_1,k_2)$ the $k_1,k_2= 1,2,3,\ldots,\infty$ full scan and, $(k_1,k_2)$ the $k_1= 1,2,3,\ldots,\infty$, $k_2= 1,2,3,\ldots,k_1$ restricted scan. We should note that interchanging $k_1$ and $k_2$ can also be defined. The small coefficients (namely the 2 and 3) in Eqs.6-7 allow easy discussion: using full scan, the divisor 3 determines all nine possible cases for $(a,b)$ pairs with $a=3k_1$, $3k_1-1$, $3k_1-2$ and $b= 3k_2$, $3k_2-1$, $3k_2-2$ by combining all *a*'s with all *b*'s. In all nine cases (see Appendix 3 for Eq.6). If *n* is integer then Eq.6 is a selector. If not, that is an irrelevant case, so we can state the following theorem:

**Theorem 3:**

All odd non-primes can be selected from the *6n+1* subset from set $S_6$ by employing a full scan if $(a,b)=(3k_1,3k_2)$ or $(a,b)=(3k_1-1,3k_2-1) \Rightarrow$ Eq.6 yields

$$n= 6k_1k_2 +k_1 +k_2 \quad \text{or} \quad 6k_1k_2 -k_1 -k_2, \qquad \text{(Eq.8)}$$

respectively, and for those *n* values, the *6n+1* is an odd non-prime; for all other *n* values which are jumped over, *6n+1* is an odd prime.

**Example for Theorem 3:**

In the case of Eqs.6 or 8 let $k_1=k_2=1$ be $\Rightarrow (a,b)=(3k_1,3k_2)=(3,3)$ yields *n*=8, and indeed 6*8+1=49 (=7*7) is non-prime, as well as that, $(a,b)=(3k_1-1,3k_2-1)=(2,2)$ yields *n*=4, and indeed 6*4+1= 25 (=5*5) is non-prime.

**Corollary 1 of Theorem 3:**

Substituting Eq.8 into *6n+1* yields:

$$6n+1= 6(6k_1k_2 +k_1 +k_2) +1 = (6k_1+1)(6k_2+1) \qquad \text{(Eq.9)}$$

and

$$6n+1= 6(6k_1k_2 -k_1 -k_2) +1 = (6k_1 -1)(6k_2 -1), \qquad \text{(Eq.10)}$$

respectively, in accordance with the fact that the product of two numbers in the form of $6n\pm1$ is also in the form of $6n\pm1$ (form conservation) and, more importantly, explains why these numbers cannot be prime (because product of two larger than 1 numbers), and correctly singled out by Eq.8. For example, the above mentioned 49 times e.g. the prime 7 yields the odd non-prime 49*7=343, but it also decomposes in the form of $S_6$ as 343= 6*57+1, as expected, so Eq.8 or 9 (with $k_1$=8, $k_2$=1) also singles out 343 from $S_6$. We also note that, in Eqs.9-10 the odd $6k_i\pm1$ can be prime or non-prime, also if $k_1=k_2$ then it is a square number $(6k_1\pm1)^2$ as demonstrated above with 49 and 25, generally, the $(6k_1\pm1)^2$ numbers are odd non-primes in the set $S_6$.

**Corollary 2 of Theorem 3:**

Both equations in Eq.8 are symmetric for the interchange of $k_1$ and $k_2$, so instead of the full scan, the restricted scan or its index-counterpart is enough for both, behaving similarly to Eq.3 in this respect, avoiding a double calculation in the scan, however, both equations are necessary to scan for all odd non-primes. For the other (see Appendix 4 for Eq.7), if *n* is an integer then Eq.7 is a selector, and if not, then it is an irrelevant case, and we can state the following, similar theorem:

**Theorem 4:**

All odd non-primes can be selected from the *6n-1* subset from the set $S_6$ with a full scan, providing $(a,b)=(3k_1,3k_2-1)$ or $(a,b)=(3k_1-1,3k_2) \Rightarrow$ Eq.7 yields

$$n= 6k_1k_2 -k_1 +k_2 \quad \text{or} \quad 6k_1k_2 +k_1 -k_2, \qquad \text{(Eq.11)}$$

respectively, and for those *n* values, *6n-1* is an odd non-prime, for all other *n* values which are jumped over, the *6n-1* is an odd prime.

Notice the strong analogy between Theorems 3 and 4 with help of Appendices 3 and 4, so the proof of Theorem 4 is analogous to Theorem 3.

**Theorem 5:**

Theorems 3-4 based on Eqs.8 and 11 select all odd primes, except 3, on account of their "jumping over" properties.

**Example for Theorem 4:**

In the case of Eq.7 let $k_1=k_2=1$ be $\Rightarrow$ the $(a,b)= (3k_1,3k_2-1)= (3,2)$ yields $n=6$, and indeed 6*6-1=35 (=5*7) is non-prime, as well as that $(a,b)= (3k_1-1,3k_2)= (2,3)$ also yields $n=6$.

**Corollary 1 of Theorem 4:**

Substituting Eq.11 for $6n-1$ yields

$$6n-1= 6(6k_1k_2 -k_1 +k_2) -1 = (6k_1+1)(6k_2-1) \qquad \text{(Eq.12)}$$

$$6n-1= 6(6k_1k_2 +k_1 -k_2) -1 = (6k_1-1)(6k_2+1) \qquad \text{(Eq.13)}$$

, respectively, in accordance with the fact that the product of two numbers in the form of $6n\pm1$ is also in the form of $6n\pm1$ (form conservation), explaining why these numbers cannot be prime, and correctly singled out by Eq.11; e.g. the above mentioned 35. We also note, that in Eqs.12-13 the odd $6k_i+1$ can be prime or non-prime, also if $k_1=k_2$ then it is true that $(6k_1+1)(6k_1-1)=36k_1^2-1$ form odd non-prime number in the set $S_6$.

**Corollary 2 of Theorem 4:**

The two equations in Eq.11 constitute dual partners in Eqs.12-13, so the indexing via a full scan generates the odd non-primes doubly for the off-diagonal elements $k_1 \neq k_2$ via interchanging values $k_1$ and $k_2$. The restricted scan removes this, but both equations are necessary in Eqs.12-13, although neither of them are symmetric with regards to interchanging $k_1$ and $k_2$, but the interchange switches one with the other and vice versa so, with respect to the pair of equations, the same pair is obtained by interchanging $k_1$ and $k_2$. The dual partnership in Eq.11 provides us with an even more important result with the scan: if a full scan is preferred over the restricted scan, then only one of the equations in Eq.11 is necessary, and the other can be discarded; in as much as interchanging $k_1$ with $k_2$ does not yield the same equation but does yield the same set of $n$ values via the scan. Finally, when using a full scan, only one of the equations in Eq.11 is necessary for the scan but, when using a restricted scan both equations in Eq.11 are necessary for the scan.

**Remarks on Theorems 1-5:**

For $S_2$ and $S_6$ we have said that $S_6$ lists the odd primes faster, (listing only, because there is no total selection) by coefficient 6 vs. 2, as well as that it jumps over the numbers in the form of $6n+3$, although $S_6$ has two equations by sign ± as opposed to only one in $S_2$. However, regarding the speed of the total selection, Eq.1 vs. Eqs.8 and 11 one should not forget that more equations, i.e., 1 vs. 3 or 4 do have to be treated - compare Tables 1 and 3 regarding this;. furthermore, the close similarity of the algebraic forms in Eqs.1, 8 and 11 is obvious, as well as in the cases of $S_2$ and Eq.1 - allowing zero for $k_1$ and/or $k_2$, Eqs.8 and 11 lose their selective property for $S_6$, and what is more, some of them generate negative, irrelevant *n* values or unselected (full) sets, i.e., $S_6$. In retrospect, starting from the right most sides of Eqs.9-10 and 12-13, the selector functions in Eqs.8 and 11 can be achieved in the same way that led us to Eq.1, but in the same way we have reached Eqs.9-10 and 12-13, it is also proven that only the forms on the right most sides have to be singled out, no other forms in the set of $6n\pm1$ exist as odd non-primes.

In view of Eqs.8-13, the twin-primes (those *n* values for which $6n+1$ and $6n-1$ are both primes) shall be commented upon at this stage but, will be discussed in detail in section 3.3. Eqs.8-13 have another important message: not only does the product of two $6n\pm1$ form odd prime or odd non-prime numbers, but also in the form of $6n\pm1$,such a trivial odd non-prime (because of the product), but the reverse is also true, such as in set $S_6$, the odd non-primes are those that can be factorized into the product of two numbers in $6n\pm1$ form; for example, the 15 is not a $6n\pm1$ form number, so none of the odd non-primes in the set $S_6$ can be divided by 15, etc., and indeed the 15 has the excluded $6n+3$ form via $n$=2, completing the total agreement in this respect.

Finally, let us take a look at the set of odd numbers in $S_2$ generated by $2k+1$ as {1,3,$\underline{5}$,$\underline{7}$, 9,$\underline{11}$,$\underline{13}$, 15,$\underline{17}$,$\underline{19}$, 21,$\underline{23}$,$\underline{25}$, 27,$\underline{29}$,$\underline{31}$,…}, wherein the $6n\pm1$ generates (underlined) all the odd primes (except 3) and certain odd non-primes, while $6n+3$ generates certain odd non-primes (not-underlined). It is therefore clear that $6n\pm1$ reduces the generation (set elements) of $2k+1$ by 1/3. For the selection, Eq.1 has only one equation in contrast to 3 or 4 in Eqs.8 and 11 (depending on whether a full or restricted scan is used), so one can say that Eq.1 is more compact and more elegant (in terms of simplicity) than Eqs.8 or 11. From this point of view the $6n\pm1$ generation of odd primes is not as "prestigious" as it would seem at first glance.

**Demonstration of Theorems 3-5**

The prime selector ability of Eqs.8 and 11 is shown in Tables 3-4 and Figure 3. Importantly, the columns containing words PRIME and TWIN in Tables 2 and 4 (or their extended ones, available on the web up to a trillion) looks very erratic in distribution, however, in view of Eqs.1-2, 8 and 11 plotted on Figures 1-3 their distribution is rather systematic. . .

### 3.2 Theorems 1-5 in relation to prime number theorem

Without going into detail, we draw attention to the way Eqs.8-13 single out odd non-primes from set $S_6$ which has an expected consequence in agreement with the density of primes stated by the well known "prime number theorem" that is, the number of primes $\leq x$, the so called $\pi(x)$ behaves asymptotically as $\sim x/ln(x)$, empirically discovered by Gauss, circa 1793 and by Legendre in 1798, which states that the density of primes in the vicinity of a natural (or real) number $x$ decreases by the value of $x$ [13-14]. Via Eqs.9-10 and 12-13 with the property that forms $6n\pm1$ but is sustained under multiplication, the odd non-primes $6n\pm1 = \Pi_{i=1..N}(6k_i\pm1)$ with $k_i$=1,2,…M, N>1 and M$\rightarrow\infty$; for order numbers n is generated more and more frequently in $S_6$ as M and N increase. If N or M is greater, $6n\pm1$ is obviously greater along $n$. The multiple choice by ± in the product Π increase the possible values for $n$, so the expression is a bit unmindful, as well as this, factors ($6k_i\pm1$) can either be primes or non-primes – better to restrict it to primes. However, one thing is obvious, the greater N and M become, the more (not only higher) $n$ values are singled out as odd non-primes in the set $S_6$, which, is in accordance with the prime number theorem, that – in this generator – the probability (or density) that $6n\pm1$ is an odd prime decreases as $n$ increases. Of course, Eq.1 behaves in the same way, in this respect, as Eqs.8 and 11, recall that form $2k+1$ is also sustained under multiplication. All in all, with Eqs.2, 9-10 and 12-13 it may be possible to change the probability approach of the density of primes to the more desired analytical approach.

Prime numbers serve as certain infinite basis sets in product generation of non-primes among natural numbers, even the basic definition of linear algebra can be recognized (see Appendix 5). Unlike in one of the most important 3 dimensional geometrical or physical spaces, for example, wherein there are three (a small finite number) unit vectors constituting the basis set, the primes as well as their density constitute an infinite basis set, with a mixed metaphor it is like the "Milky Way": to get an idea, in the interval [1,500] there are 95 primes, the smallest neighboring difference is 3-2=1, the largest is 127-113=14, but further away, (of course "far" is relative) with the same size interval [$10^{12}$-500, $10^{12}$] there are still 21 primes. The smallest neighboring difference is between twin primes 999999999961 - 999999999959 = 2, the largest is 999999999847 - 999999999847 = 80, (keep in mind that by the hypothesis that there are infinite twin primes, the smallest difference 2 can, in principle, be found up to infinity so, the largest difference refers to the decrease of density of primes when considering a fixed finite interval). Although these last properties, among many others, are well known, we mention this, because these properties visibly belong to the complement image of selector functions in Eqs.1-3 and are commented upon in section 2.2 or with Eqs.8 and 11.

### 3.3 Twin primes with selector functions

Eqs.1-2 provide the selector function for set $S_2$ for the odd number generator $2k+1$, and twin primes are found if a $k_2$ is found that Eq.1 jumps over, along with jumping over $k_2+1$, so $2k_2+1$ and $2(k_2+1)+1= 2k_2+3$ are the so called twin primes with this special $k_2$ found by scanning $a$ and $b$ (in the image of Eq.1), - index 2 indicates "twin". Important theorems are: 1. There are infinite (odd) primes (mentioned above). 2. There are fewer twin primes than isolated primes (trivial, not because one twin prime contains two primes, but because there are isolated primes, e.g., 23) and, 3. The hypothesis is that there are infinite twin primes (not proved yet in the literature). In $S_6$ the $6n\pm1$ generates all the odd primes, but unfortunately some odd non-primes too. Another property is:

**Theorem 6:**

No triplet primes (neighboring primes of three) except the triplet {3,5,7} or higher multiples, – actually, beside the existence of an infinite set of twin primes this exception is a quintet: {2,3,5,7}.

**Proof of Theorem 6:**

It can be proved easier with the $6n\pm1$ generator than with Eq.1 and it follows immediately, because there is an intermediate odd non-prime number between any order number $n$ and $n+1$, that is in the above mentioned interval $[6n+1, 6(n+1)-1]= [6n+1, 6n+5]$, the odd non-prime $6n+3$, so if $6n+1$ and $6n+5$ happen to be two odd primes in accordance with the selection property of Eqs.8 and 11, the $6n+3=3*(2n+1)$ odd non-prime always cancels the opportunity for higher multiplicity than two (twin), - in fact the $6n+3$ is a selector again, as it was once before in Lemma 1, it singles out (or cancels the opportunity of being) a higher multiplicity than the twin primes above number/prime 7.

**Theorem 7:**

If both $6n_2\pm1$ are odd primes for an $n_2$, those are the very twin primes we are discussing, and Eqs.8 and 11 are the selector functions which must jump over this isolated $n_2$ via scanning $(k_1,k_2)$ for the four equations in Eqs. 8 and 11 simultaneously.

**Extending Theorem 7:**

In view of the twin hypothesis, these $n_2$ values in Theorem 7 constitute the subset only of natural numbers and have no upper boundary, i.e., $n_2\to\infty$, - index 2 also indicates "twin". On the other side, there must be an infinite number of, note them as, $n_0$ values for which Eq.8 and 11 both hit this value), so both $6n_0\pm1$ are odd non-primes. Because the density of odd primes decreases as shown by the prime number theorem, more and more of these $n_0$ exist among the larger $n$ values, - index zero indicates that none of $6n_0\pm1$ is prime. In conclusion, there are, (note it as, $n_1$ values for which one (indicated in the index) of $6n_1\pm1$ is an odd prime and the other is not in the pair neighborhood) a finer subdivision comes from that of the upper or lower which is the odd prime. Besides the fact that twin primes i.e., twin odd primes are in the focus of research, (noted here in the index in $n_2$) the distribution and appearance of indices just defined, $n_1$ (isolated odd primes, i.e. non-twin primes) and $n_0$, are also interesting, however this is not discussed any further here.

**Demonstration of Theorem 7:**

Finally, with language of logic for the evaluation among integer numbers: 1. If for a given $n$ a not necessarily unique $(k_1,k_2)$ value is found for any of the two equations in Eq.8 to hold $\Rightarrow$ $6n+1$ is odd non-prime; and for any of the two equations in Eq.11 $\Rightarrow$ the $6n-1$ is an odd non-prime. 2. If no $(k_1,k_2)$ exists to yield this $n$, i.e., that the scan jumps over via any of Eq.8, then $6n+1$ is odd prime; via any of Eq.11, then $6n-1$ is odd prime. 3. If no $(k_1,k_2)$ exists to yield this $n$ i.e., that the scan jumps over via any of the four equations in Eqs.8 and 11, then $6n\pm1$ is odd twin prime at this $n$. Table 4 and Figure 3 demonstrate this selector ability for twin primes.

A note on Theorem 7 is that, obviously $(6n_2+1)/(6n_2-1)\to1$ if $n_2\to\infty$ irrespectively of the property of $n_2$ above. However, it would be interesting to figure out if the techniques here can be benefited to study the ratio of any consecutive prime numbers [15] (not only the subset twin primes), it will be discussed in another study.

### 4.1 Generalization

Any natural number can be decomposed with primes such as $n= 2^{n1}3^{n2}5^{n3}...$ with $n_i \geq 0$ natural numbers (Appendix 5), where the first term provides the odd ($n_1$=0) or even ($n_1$>0) property, and the rest is always odd (including 1 via $n_i$=0 for all i>1) yielding form $n=2^{n1}m_{odd}$. The generalization for the previous sections above comes from triple or the more general *n*-tuple products $m= (2a+1)(2b+1)(2c+1)$ or $\Pi_{i=1...N} = (2a_i+1)$ for $m_{odd}$, yielding the polynomial series

$$m_1 = 2a+1 \quad \text{(Eq.14.1)}$$

$$m_2 = (2a+1)(2b+1)= 4ab +2(a+b) + 1 \quad \text{(Eq.14.2)}$$

$$m_3 = (2a+1)(2b+1)(2c+1) = 8abc +4(ab+ac+bc) +2(a+b+c) +1 \quad \text{(Eq.14.3)}$$

$$m_N = \Pi_{i=1...N}(2a_i+1)= 2^N(\Pi_{i=1...N}\, a_i) +... +2(\Sigma_{i=1...N}\, a_i) +1, \quad \text{(Eq.14.4)}$$

where indices on *m* show the number of terms in the product and provide the following definition.

**Definition 1:**

Let the generalization of the selector function *k* in Eq.1 be $2ab +a+b$, $4abc +2(ab+ac+bc) +a+b+c,...,$ $2^{N-1}(\Pi_{i=1...N}\, a_i) +... +(\Sigma_{i=1...N}\, a_i)$, respectively, adopted from Eq.14.

What is important is that the terms in the products in Eq.14 can be odd primes or composite odd numbers, and one should bear in mind again, that the product of numbers in the form of $(2a_i+1)$ preserves this form. The parameters run as $a, b, c, a_i$ =1,2,3,… ∞, otherwise, e.g., *c*=0 makes Eq.14.3 to fall into Eq.14.2, etc.; compare to full scan or generalize the restricted scan. How do these symmetric (with respect to interchange of any $a_i$ and $a_j$) polynomials select natural numbers or associate to the number of their divisors? The $m_1$ lists all the odd numbers (odd primes and odd non-primes), $m_2$ (the same as Eq.2 extended with index just defined) lists all odd non-primes, i.e., odd numbers which can be decomposed to a product of two or more (non-unit) numbers, or jumps over numbers (odd primes) which cannot be decomposed to a product of two or more (non-unit) numbers. For example, the starting values *a*=*b*=*c*=1 and *a*=*b*=1 with *c*=2, yield $m_3$= 27= 9*3= 3*3*3 and 45= 9*5= 3*3*5, respectively, so, these can be decomposed at least three terms, and jump (from the beginning) over 3, 5, 7, 9=3*3, 11, 13, 15=3*5, 17, 19, 21=3*7, 23, 25=5*5, (27 is hit), 29, 31, 33=3*11, 35=5*7, 37, 39=3*13, 41, 43, (45 is hit), i.e.,those numbers which cannot be decomposed over three or more non-unit terms, but maximum one- (e.g., the 3, 9) or two (e.g. the 9). We have arrived at the generalization of Theorems 1-5:

**Theorem 8:**

The polynomial $m_N$ in Eq.14 jumps over those odd numbers which can be decomposed only to maximum 1 or 2,…, or N-1 non-unit terms, but cannot be decomposed to a product of N or more non-unit terms and, hits those odd numbers which can be decomposed to product N or more non-unit terms. (Of course, the self is not considered divisor or decomposition.)

**Corollary of Theorem 8:**

In this way the term or definition of odd primes can be generalized/extended, selected and characterized by Eq.14: "1st, 2nd, 3rd, …Nth order" odd primes are those which only have a maximum of: 1,2,3,…N, and no more non-unit divisors, respectively (neglecting the 1 and the self as divisor) the 3, 5, and 7 are 1st order odd primes (conventional primes); the 9 and 15 are 2nd order odd primes, etc, of course the most prestigious is the 1st order prime. The conventional odd primes, and Eq.14 provides information on their exact location and distribution, at least it converts or associates the number of divisors of natural numbers to polynomials in Eq.14.

The generator $4n\pm1$ has been mentioned in section 3.1. In relation to its generalization in this section, the following can be stated: with $k$, $n$= 1,2,3,… formulas $2k+1$ and $4n\pm1$ generate all odd natural numbers, while $6n\pm1$ generates all odd natural numbers except those which can be divided by 3, more importantly, the odd non-primes are selected by Eq.1 for $2k+1$, Eqs.8 and 11 for $6n\pm1$, leaving the residue set containing all the odd primes and only those, discussed in detail above. For the $4n\pm1$ rule we can analogously state:

**Theorem 9:**

In a likewise fashion to Theorems 3-5, $n= 4k_1k_2 +k_1 +k_2$ and $4k_1k_2 -k_1 -k_2$ selects all odd non-primes generated by $4n+1$; $n= 4k_1k_2 +k_1 -k_2$ and $4k_1k_2 -k_1 +k_2$ selects all odd non-primes generated by $4n-1$, and the remainder set generated by $4n\pm1$ contains all the odd primes and only those via a full scan, as well as the properties (or, proof) listed for $6n\pm1$ above or in Appendices 3-4) this can be analogously stated for $4n\pm1$, too.

This work was fully presented [footnote 8] in 14th International Conference of Numerical Analysis and Applied Mathematics (ICNAAM 2016).

**Conclusion**

It is not the prime number generator, but the counterpart generator functions that have been arrived at, introduced, and analyzed - as indicated in the Abstract - which completely single out the odd non-primes from the set of odd numbers having at least, simple and analytic forms [16]. The author looks at these selector functions and sometimes deems them quite trivial, and sometimes extremely powerful, but one thing is certain, these items cannot be neglected in the ongoing discussion on prime number generating functions. There have been many tables of primes produced over the years, and it is said that they show that the detailed distribution is quite "erratic", but the number of primes up to positive real or integer $x$ in a residual set (all odd primes) that the $m= 2(2ab +a +b)+1$ selector function leaves from the set of all odd numbers (plus 1 for the prime 2), the so called $\pi(x)$ grows fairly steadily. In view of the existing, discovered and introduced selector functions in the Abstract, which systematically single out all the odd non-primes from odd natural numbers leaving all odd primes, more, or by request only the twin primes in the set, the word "erratic" should be revised.

As $x^2+1=0$ generates one of the most characteristic number (imaginary unit, $i\equiv\sqrt{-1}$) of algebra via its enforced solution, $k=2ab+a+b$ generates (jumps over) the most characteristic numbers of number theory, the order number of odd primes via scanning $(a,b)$.

**Appendices**

Appendix 1.: A scan for an individual $k_0$ in Eq.1 with restricted scan yields the straightforward interval as follows: the minimal $a=M_1$ gives maximal value if $a=b=M_1$ so, according to Eq.1 the $k_0\approx 2M_1^2+2M_1$, where the approximate sign indicates that there may not be an integer solution; it yields $0<M_1=((1+2k_0)^{1/2}-1)/2$, which is $\approx\sqrt{(k_0/2)}$ for large $k_0$. The $\{M_1\}$ is the smaller neighbor integer, e.g., $\{3.14\}=3$, and "a" must start with $\{M_1\}$. The maximal $a=M_2$ value gives minimal value if $b=1$, so, going by Eq.1 $k_0\approx 3M_2+1$, yielding $M_2=(k_0-1)/3$ - which is $\approx k_0/3$ for large $k_0$. From these the range for "$a$" is $\{M_1\}\le a<[M_2]$, where $[M_1]$ is the larger neighbor integer, e.g., $[3.14]=4$. For $a=\{M_1\}$ the $b=1,2,\ldots,\{M_1\}-1$ gives $k<k_0$, so $b=\{M_1\}$ should be calculated only, for $a=[M_2]$ and $b=1$ gives $k>k_0$, so "$a$" must terminate with $[M_2]-1$. For intermediate $a$ values $2ab+a+b\approx k_0$ yields $b\approx(k_0-a)/(2a+1)$, so with integer values $b$ runs for the $\{(k_0-a)/(2a+1)\}$ and $[(k_0-a)/(2a+1)]$ values yielding neighboring under or over values, or may hit the $k_0$ itself. For an interval for $k$, it can be derived in the same way, avoiding a double calculation for neighboring $k$ values. (The expressions for $M_1$ and $M_2$ are in accordance with the initial range a natural number can be searched for its possible divisors.); e.g., for $k_0=20000$, the $\{M_1\}=99$ and $[M_2]=6667$, for $a=99$, the $b=100$ and 101, yielding $k=19999$ and 20198 i.e., under and over 20000 via Eq.1, as well as for $a=6666$, the $b=1$ and 2, yielding $k=19999$ and 33332, i.e., under and over 20000. During the scan $a=1538$ require $b=6=\{6\}=[6]$ integer, for which $k=k_0=20000$, and indeed $2k_0+1=40001=13*3077$ is not an odd prime.

In other words, $m=2k+1$ is the $k^{th}$ odd natural number, if exits integer $(a,b)$ such as $k=2ab+a+b$ for any (one or more) of $a,b=1,2,\ldots$ value pairs, then $m$ is an odd non-prime. Since $b=(k-a)/(2a+1)\ge 1$, then $(k-1)/3\ge a$, so for $a=1,2,\ldots,((k-1)/3)$ if all $b=(k-a)/(2a+1)$ is non-integer, then $m$ is odd prime, if at least one of it is integer, then $m$ is odd non-prime. For example, $m=33$ $\Rightarrow k=16$, for $a=1,2,\ldots,((k-1)/3)=(15/3)=5$, then $b=5$, 14/5, 13/7, 12/9, 1 $\Rightarrow$ there exist integer in this set, the 1 and 5, so $m=33$ is an odd non-prime, and indeed $m=33=3\text{x}11$. Another example, $m=31\Rightarrow k=15$, for $a=1,2,\ldots,((k-1)/3)=(14/3)=4$, then $b=14/3$, 13/5, 12/7, 11/9 $\Rightarrow$ there is no integer in this set, so $m=31$ is an odd prime, and indeed it is.

Appendix 2.: One should keep in mind that although Eqs.1-2 are two dimensional polynomials, but it does not mean that it is possible to describe the complement of image $m$ (odd primes) fully with two or higher dimensional polynomials. For example, Fermat numbers with $q=0$-4 is described with powers and not polynomials (although far not effectively, because the very few values for $q$), however contrary, the Euler polynomial exists and more powerful in that it generates a few more primes, but far not all the odd primes. However, just keeping the hope, powers of odd non-primes can be generated with polynomials in Eq.2 easily, e.g. $5^2=m(a,b)$ is satisfied with $a=b=2$, see Table 1. Contra reasoning again, one can call the textbook example of $\int(dx/\ln(x))$ of which primitive function exists (like odd primes exist) but no finite expression of analytical functions describes it – i.e. the dilemma such as generator function for odd primes exists or not is still an issue, as well as Eq.2 is yet "mute" about the form for describing the complement of Im($m$).

Another example in this respect here is didactic or rhetoric: For example, the $\kappa\equiv 2^a$ selects "power of 2" with $a$=1,2,3,…, and it has nothing to do directly with the order number of odd primes (although if "$a$" is prime, the number before it is a Mersenne prime, for example), but selects from natural numbers as Eq.1 does from odd numbers with simple function. However, to describe the complement of image set Im($\kappa$) is not simple function, although it is simple to list consecutively: {$a$= 1,2,3,4,5,6,7,8,9…}\{$2^a$= 2,4,8,…} = {1,3,5,6,7,9,…}. And the line of examples along this way can be continued easily, e.g. with the above mentioned diagonal elements $m$ in Table 1: if $\kappa := m = (2a+1)^2$, what is the generator function for complement of Im($\kappa$)?

Appendix 3.: For Eq.6, if (a,b)=($3k_1$,$3k_2$) $\Rightarrow$ n=(2ab+a+b)/3 = ($2*3k_1*3k_2$ +$3k_1$ +$3k_2$)/3= $6k_1k_2$ +$k_1$ +$k_2$= integer $\Rightarrow$ selector, if (a,b)= ($3k_1$,$3k_2$-1) $\Rightarrow$ n=($2*3k_1*(3k_2-1)+3k_1+(3k_2-1)$)/3= $2k_1*(3k_2-1)+k_1+k_2-1/3$= non-integer, if (a,b)= ($3k_1$, $3k_2$-2) $\Rightarrow$ n=integer-2/3 = non-integer, if (a,b)= ($3k_1$-1,$3k_2$ ) $\Rightarrow$ same as (a,b)= ($3k_1$,$3k_2$-1) $\Rightarrow$ n= non-integer, if (a,b)= ($3k_1$-1,$3k_2$-1) $\Rightarrow$ n= $6k_1k_2$ -$k_1$ -$k_2$= integer $\Rightarrow$ selector, if (a,b)= ($3k_1$-1,$3k_2$-2) $\Rightarrow$ n= integer +1/3= non-integer, if (a,b)= ($3k_1$-2,$3k_2$ ) $\Rightarrow$ same as (a,b)= ($3k_1$,$3k_2$-2) $\Rightarrow$ n= non-integer, if (a,b)= ($3k_1$-2,$3k_2$-1) $\Rightarrow$ same as (a,b)= ($3k_1$-1,$3k_2$-2) $\Rightarrow$ n= non-integer, if (a,b)= ($3k_1$-2,$3k_2$-2) $\Rightarrow$ n= integer + 4/3 $\Rightarrow$ n= non-integer.

Appendix 4.: For Eq.7, if (a,b)=($3k_1$,$3k_2$) $\Rightarrow$ n=(2ab+a+b+1)/3 = ($2*3k_1*3k_2$ +$3k_1$ +$3k_2$+1)/3= integer+1/3 = non-integer, if (a,b)= ($3k_1$,$3k_2$-1) $\Rightarrow$ n=($2*3k_1*(3k_2-1)+3k_1+(3k_2-1)+1$)/3= $6k_1k_2-k_1+k_2$= integer $\Rightarrow$ selector, if (a,b)= ($3k_1$, $3k_2$-2) $\Rightarrow$ n=integer-1/3 = non-integer, if (a,b)= ($3k_1$-1,$3k_2$ ) $\Rightarrow$ same as (a,b)= ($3k_1$,$3k_2$-1) $\Rightarrow$ n= $6k_1k_2+k_1-k_2$= integer $\Rightarrow$ selector, if (a,b)= ($3k_1$-1,$3k_2$-1) $\Rightarrow$ n= integer +1/3= non-integer, if (a,b)= ($3k_1$-1,$3k_2$-2) $\Rightarrow$ n= integer +2/3= non-integer, if (a,b)= ($3k_1$-2,$3k_2$ ) $\Rightarrow$ same as (a,b)= ($3k_1$,$3k_2$-2) $\Rightarrow$ n= non-integer, if (a,b)= ($3k_1$-2,$3k_2$-1) $\Rightarrow$ same as (a,b)= ($3k_1$-1,$3k_2$-2) $\Rightarrow$ n= non-integer, if (a,b)= ($3k_1$-2,$3k_2$-2) $\Rightarrow$ n= integer + 5/3 $\Rightarrow$ n= non-integer.

Appendix 5.: Not going to detail et all, only an example is given for a connection between primes and linear algebra: By the "fundamental theorem of arithmetic", with the set of primes {2,3,5,7,11,…} a coordinate ($n_1$,$n_2$,…) can be unequivocally defined for any natural number $n$= $2^{n1}3^{n2}$… where $n_i$ are also natural numbers including zeros, e.g. 63= $3^2$*7= (0,2,0,1,0…), and $63^2$= $3^4$*$7^2$= 2*(0,2,0,1,0…)= (0,4,0,2,0…) indeed, or 5*63= (0,0,1,0,…) + (0,2,0,1,0…) = (0,2,1,1,0…) indeed, - a connection between products among natural numbers as well as sum and constant product rule with infinite Euclidean basis set. Notice that the "fundamental theorem of arithmetic" mentioned in the beginning in section 2.1 is only recomposed and extended here with the inclusion of zeros for $n_i$'s.

**Acknowledgement**

Support and inspiration from grant OTKA--K 112312, 115733 and 119358 is warmly appreciated.

**Footnotes:**

1.: http://en.wikipedia.org/wiki/Prime_number

Figures 1 legend: As example, more extended interval [1,100] than listed in Table 2 is plotted exhibiting the selector function *k(a,b)= 2ab+a+b* for selecting all odd primes via K=0 values (those of which *k(a,b)* is hit zero times). Notice the total coincidence such as K=0 hits all odd primes (all primes except 2) in this "spectrum", in a certain view it is much more powerful then the Euler prime generator which hits only 40 almost consecutive odd primes in the beginning of the set of natural numbers.

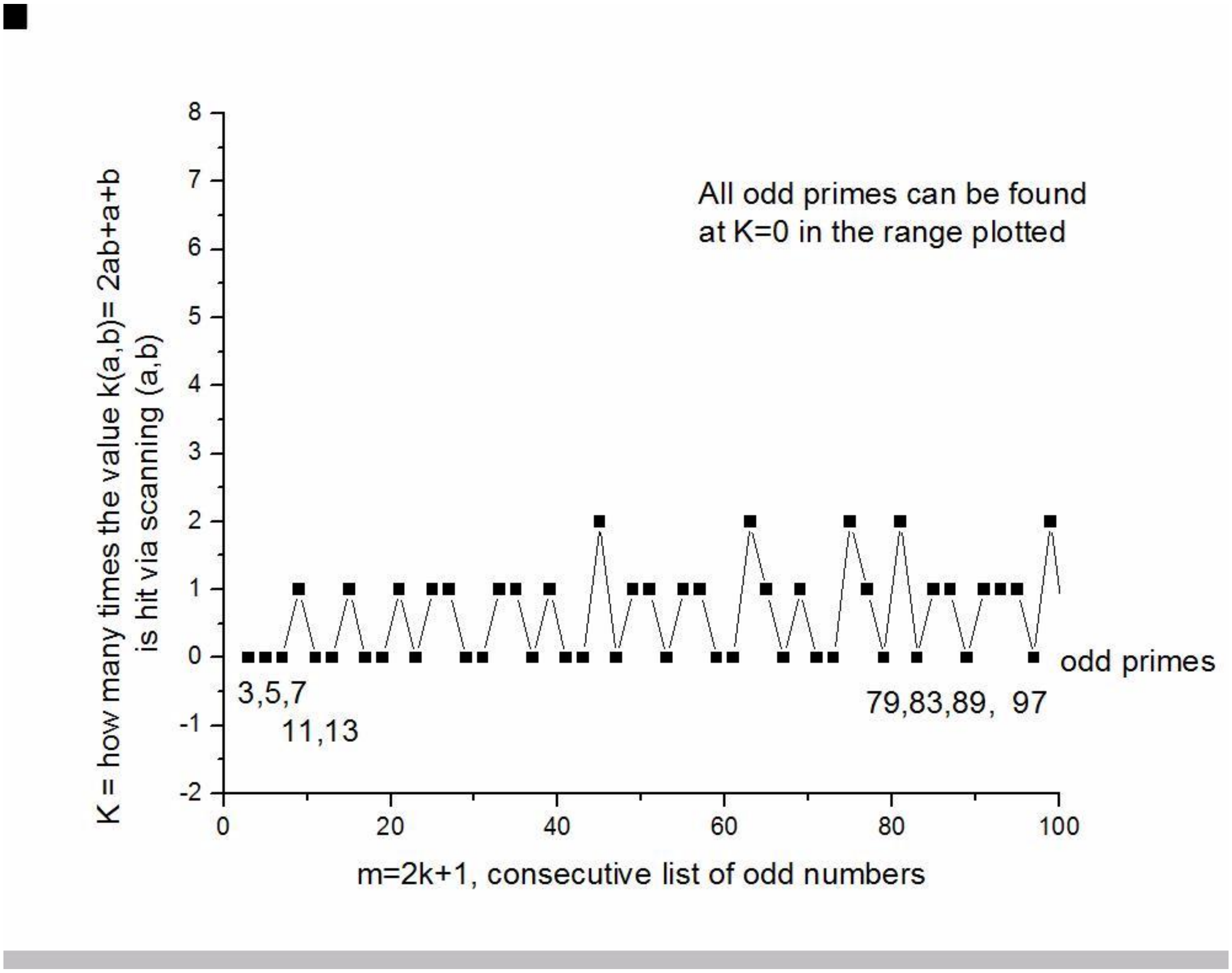

Figure 2 legend: Same as Figure 1, but for interval [901,1000].

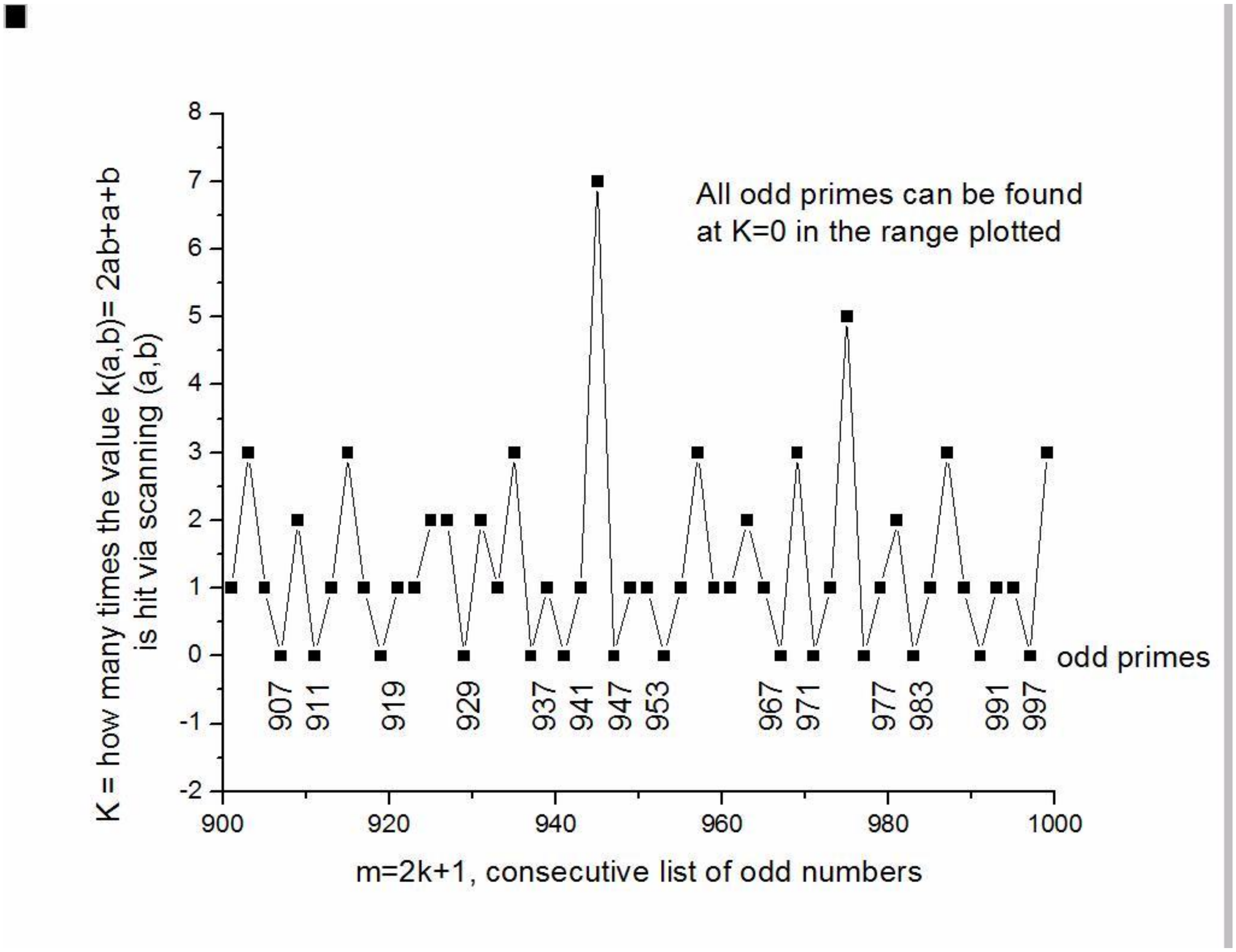

Figure 3 legend: Extended interval than listed in Table 4 is plotted exhibiting the power of selector functions in Eqs.8 and 11 for selecting twin primes via $K_-+K_+=0$ values; see Table 4 caption for the definition of $K_-+K_+$. Notice the total coincidence such as $K_-+K_+=0$ hits all odd twin primes (except primes 2 and 3 in this respect) in this “spectrum”.

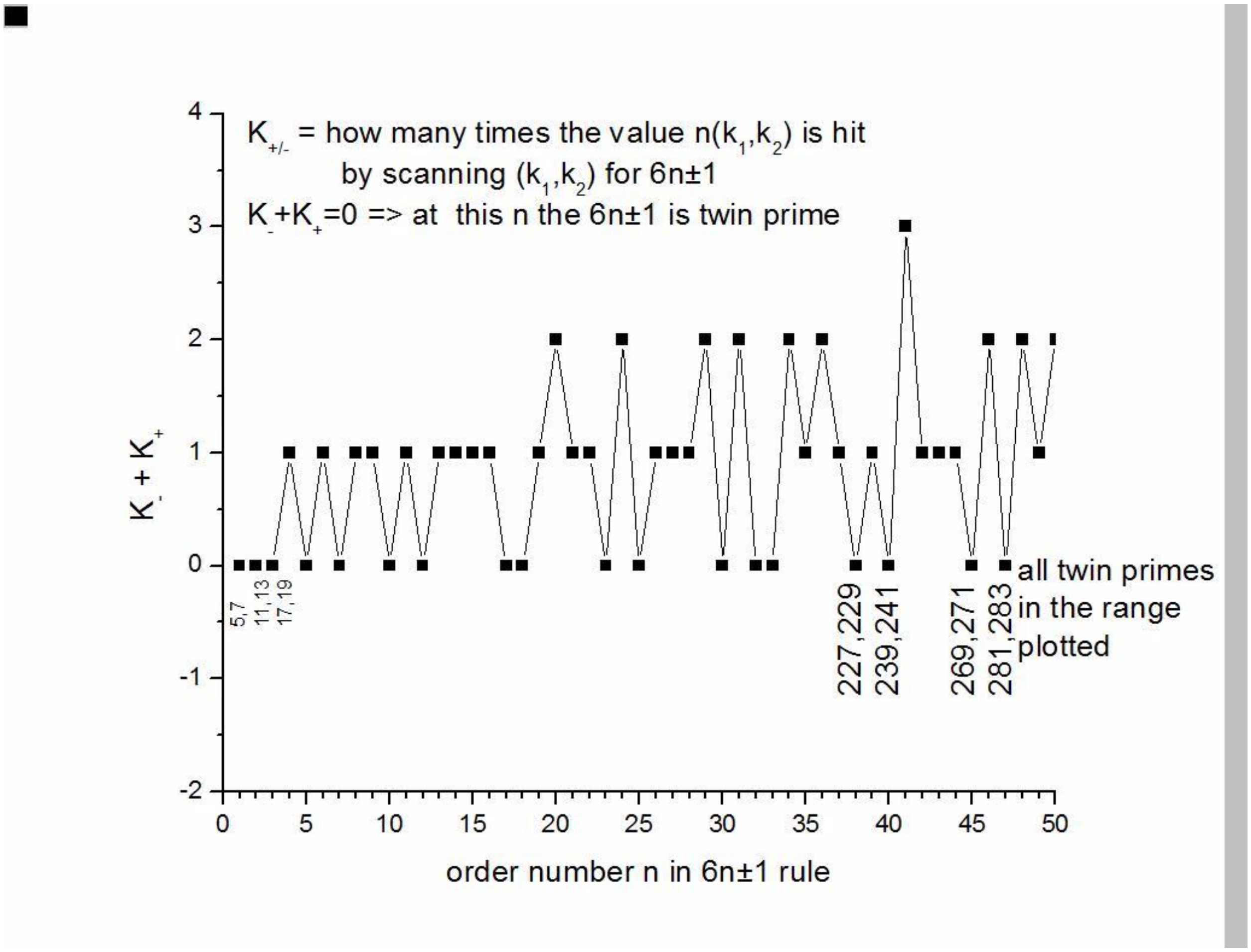

Table 1.: Beginning values of the order number or selector function *k(a,b)= 2ab+a+b* for all odd numbers (*m=2k+1*) which are odd non-primes by Eq.1 as total selection via scanning with (*a*,*b*), e.g. *a*=*b*=1 yield *k*=4 in the table, the 4th odd number, i.e. the *m*= 2*4+1= 9=3*3 odd non-prime.

*k(a,b)= 2ab+a+b* values (order numbers of odd numbers being non-primes):

| | b=1 | 2 | 3 | 4 | 5 | 6 | 7 | 8 | 9 | 10 |
|---|---|---|---|---|---|---|---|---|---|---|
| a=1 | 4 | | | | | | | | | |
| 2 | 7 | 12 | | | | | | | | |
| 3 | 10 | 17 | 24 | | | | | | | |
| 4 | 13 | 22 | 31 | 40 | | | | | | |
| 5 | 16 | 27 | 38 | 49 | 60 | | | | | |
| 6 | 19 | 32 | 45 | 58 | 71 | 84 | | | | |
| 7 | 22 | 37 | 52 | 67 | 82 | 97 | 112 | | | |
| 8 | 25 | 42 | 59 | 76 | 93 | 110 | 127 | 144 | | |
| 9 | 28 | 47 | 66 | 85 | 104 | 123 | 142 | 161 | 180 | |
| 10 | 31 | 52 | 73 | 94 | 115 | 136 | 157 | 178 | 199 | 220 |
| a | 1+3a | 2+5a | 3+7a | 4+9a | 5+11a | 6+13a | . . . | | | |

*m=2(2ab+a+b)+1* values (odd numbers being non-primes, the odd primes are jumped over):

| | b=1 | 2 | 3 | 4 | 5 | 6 | 7 | 8 | 9 | 10 |
|---|---|---|---|---|---|---|---|---|---|---|
| a=1 | 9 | | | | | | | | | |
| 2 | 15 | 25 | | | | | | | | |
| 3 | 21 | 35 | 49 | | | | | | | |
| 4 | 27 | 45 | 63 | 81 | | | | | | |
| 5 | 33 | 55 | 77 | 99 | 121 | | | | | |
| 6 | 39 | 65 | 91 | 117 | 143 | 169 | | | | |
| 7 | 45 | 75 | 105 | 135 | 165 | 195 | 225 | | | |
| 8 | 51 | 85 | 119 | 153 | 187 | 221 | 255 | 289 | | |
| 9 | 57 | 95 | 133 | 171 | 209 | 247 | 285 | 323 | 361 | |
| 10 | 63 | 105 | 147 | 189 | 231 | 273 | 315 | 357 | 399 | 441 |

Table 2.: Exhibition of the multiplicity of selector function $k(a,b)= 2ab+a+b$ (totally selecting odd non-primes from all odd numbers) in Eq.1 for beginning values; odd primes are those of which $k$ is hit zero times (jumped over), if $k$ is hit at least once, those are odd non-primes. (E.g. $k$=22 in Table 1 is hit 2 times, no more hit, since all columns in Table 1 increase downwards.) K ≡ how many times the value $k(a,b)$ is hit by Eq.1 via restricted scan; $m \equiv 2k+1$ = consecutive list of odd numbers. Larger intervals are plotted on Figures 1-2. Notice the total coincidence in columns K at K=0 at which $m=2k+1$ is odd prime (PRIME).

| k(a,b) | K | m=2k+1 | m |
|---|---|---|---|
| 1 | 0 | 3 | PRIME |
| 2 | 0 | 5 | PRIME |
| 3 | 0 | 7 | PRIME |
| 4 | 1 | 9 | |
| 5 | 0 | 11 | PRIME |
| 6 | 0 | 13 | PRIME |
| 7 | 1 | 15 | |
| 8 | 0 | 17 | PRIME |
| 9 | 0 | 19 | PRIME |
| 10 | 1 | 21 | |
| 11 | 0 | 23 | PRIME |
| 12 | 1 | 25 | |
| 13 | 1 | 27 | |
| 14 | 0 | 29 | PRIME |
| 15 | 0 | 31 | PRIME |
| 16 | 1 | 33 | |
| 17 | 1 | 35 | |

Table 3.: Beginning values (scanning with $(k_1,k_2)$) of the order number by selector function $n(k_1,k_2)$ for odd numbers in form $6n\pm1$ which are odd non-primes by Eqs.8 and 11 as total selection, e.g. $(k_1,k_2)$=(1,1) yield $n= 6k_1k_2 +k_1 +k_2$ = 8 for $6n+1$=6*8+1=49= 7*7 in the table, an odd non-prime in the form $6n+1$. On the other side, for example for $6n+1$ form, let say in the beginning range $n$=1-9, the 4, 8 and 9 are selected out, the residue is $n$=1,2,3,5,6,7, and indeed $6n+1$= 7,13,19,31,37,43, respectively are indeed $6n+1$ form odd primes, while for $n$= 4,8,9 the $6n+1$= 25,49,55, respectively are indeed $6n+1$ form odd non-primes – completing a total division for odd prime/non-prime in this range. In other words, in the beginning range listed ($0< k_1,k_2 \leq 5$), all the odd non-primes start to be selected from odd numbers in form $6n\pm1$, those odd numbers in form $6n\pm1$ are jumped over which are primes. Restricted scan is used for $n_a$ and $n_b$ as well as full scan for $n_c$, where the a,b,c indexing distinguishes for variable $n$ in this respect.

| $k_1$ | $k_2$ | $n_a=6k_1k_2+k_1+k_2$ | $6n_a+1$ | $n_b=6k_1k_2-k_1-k_2$ | $6n_b+1$ |
|---|---|---|---|---|---|
| 1 | 1 | 8 | 49 | 4 | 25 |
| 2 | 1 | 15 | 91 | 9 | 55 |
| 2 | 2 | 28 | 169 | 20 | 121 |
| 3 | 1 | 22 | 133 | 14 | 85 |
| 3 | 2 | 41 | 247 | 31 | 187 |
| 3 | 3 | 60 | 361 | 48 | 289 |
| 4 | 1 | 29 | 175 | 19 | 115 |
| 4 | 2 | 54 | 325 | 42 | 253 |
| 4 | 3 | 79 | 475 | 65 | 391 |
| 4 | 4 | 104 | 625 | 88 | 529 |
| 5 | 1 | 36 | 217 | 24 | 145 |
| 5 | 2 | 67 | 403 | 53 | 319 |
| 5 | 3 | 98 | 589 | 82 | 493 |
| 5 | 4 | 129 | 775 | 111 | 667 |
| 5 | 5 | 160 | 961 | 140 | 841 |

| $k_1$ | $k_2$ | $n_c=6k_1k_2+k_1-k_2$ | $6n_c-1$ |
|---|---|---|---|
| 1 | 1 | 6 | 35 |
| 1 | 2 | 11 | 65 |
| 1 | 3 | 16 | 95 |
| 1 | 4 | 21 | 125 |
| 1 | 5 | 26 | 155 |
| 2 | 1 | 13 | 77 |
| 2 | 2 | 24 | 143 |
| 2 | 3 | 35 | 209 |
| 2 | 4 | 46 | 275 |
| 2 | 5 | 57 | 341 |
| 3 | 1 | 20 | 119 |
| 3 | 2 | 37 | 221 |
| 3 | 3 | 54 | 323 |
| 3 | 4 | 71 | 425 |
| 3 | 5 | 88 | 527 |
| 4 | 1 | 27 | 161 |
| 4 | 2 | 50 | 299 |
| 4 | 3 | 73 | 437 |
| 4 | 4 | 96 | 575 |
| 4 | 5 | 119 | 713 |
| 5 | 1 | 34 | 203 |
| 5 | 2 | 63 | 377 |
| 5 | 3 | 92 | 551 |
| 5 | 4 | 121 | 725 |
| 5 | 5 | 150 | 899 |

Table 4.: Exhibition of the multiplicity of selector function $n(k_1,k_2)$ (totally selecting odd non-primes from generator $6n\pm1$) in Eqs.8 and 11 for beginning values; odd primes are those of which $n$ is hit zero times, if $n$ is hit at least once, those are odd non-primes. In detail: The $n$ is a consecutive list of order (natural) numbers.

$K_-$ ≡ full scanning, how many times the value $n$ is hit by (any but only one choice of) Eq.11 for $6n-1$, if $K_-=0 \Rightarrow 6n-1$ is odd prime.

$K_+$ ≡ restricted scanning, how many times the value $n$ is hit by (any but both are to be checked in) Eq.8 for $6n+1$, if $K_+=0 \Rightarrow 6n+1$ is odd prime.

As a consequence, $K_-+K_+=0 \Rightarrow$ no Eq.8 nor Eq.11 hits value $n$ in scanning $(k_1,k_2)$ via full and restricted scan $\Rightarrow 6n\pm1$ is twin prime (TWIN) at this $n$.

| n | $K_-$ | 6n-1 | $K_+$ | 6n+1 | $K_-+K_+$ | |
|---|---|---|---|---|---|---|
| 1 | 0 | 5 | 0 | 7 | 0 | TWIN |
| 2 | 0 | 11 | 0 | 13 | 0 | TWIN |
| 3 | 0 | 17 | 0 | 19 | 0 | TWIN |
| 4 | 0 | 23 | 1 | 25 | 1 | |
| 5 | 0 | 29 | 0 | 31 | 0 | TWIN |
| 6 | 1 | 35 | 0 | 37 | 1 | |
| 7 | 0 | 41 | 0 | 43 | 0 | TWIN |
| 8 | 0 | 47 | 1 | 49 | 1 | |
| 9 | 0 | 53 | 1 | 55 | 1 | |
| 10 | 0 | 59 | 0 | 61 | 0 | TWIN |
| 11 | 1 | 65 | 0 | 67 | 1 | |
| 12 | 0 | 71 | 0 | 73 | 0 | TWIN |
| 13 | 1 | 77 | 0 | 79 | 1 | |
| 14 | 0 | 83 | 1 | 85 | 1 | |
| 15 | 0 | 89 | 1 | 91 | 1 | |
| 16 | 1 | 95 | 0 | 97 | 1 | |
| 17 | 0 | 101 | 0 | 103 | 0 | TWIN |
| 18 | 0 | 107 | 0 | 109 | 0 | TWIN |
| 19 | 0 | 113 | 1 | 115 | 1 | |
| 20 | 1 | 119 | 1 | 121 | 2 | |
| 21 | 1 | 125 | 0 | 127 | 1 | |
| 22 | 0 | 131 | 1 | 133 | 1 | |
| 23 | 0 | 137 | 0 | 139 | 0 | TWIN |
| 24 | 1 | 143 | 1 | 145 | 2 | |
| 25 | 0 | 149 | 0 | 151 | 0 | TWIN |
| 26 | 1 | 155 | 0 | 157 | 1 | |
| 27 | 1 | 161 | 0 | 163 | 1 | |
| 28 | 0 | 167 | 1 | 169 | 1 | |
| 29 | 0 | 173 | 2 | 175 | 2 | |
| 30 | 0 | 179 | 0 | 181 | 0 | TWIN |